\documentclass{ifacconf}

\usepackage{graphicx}      
\usepackage{natbib}        
\usepackage{float}
\usepackage{amsmath,amssymb,amsfonts,amstext}
\usepackage{epsfig}
\usepackage{setspace}
\usepackage{hhline}
\usepackage{color,soul}
\usepackage[dvipsnames]{xcolor}
\definecolor{Light}{gray}{.90}
\usepackage[justification=centering]{caption}
\usepackage{multirow}
\usepackage{booktabs}
\usepackage{bm}
\usepackage[english]{babel}
\usepackage{tikz}
\usepackage{tikzscale}
\usepackage{epstopdf}
\usepackage{adjustbox}
\usepackage{longtable}
\usepackage{accents}
\usepackage{lineno}
\usepackage{multirow}
\usepackage{colortbl}
\usepackage{ragged2e}
\usepackage{amsmath,amssymb}
\makeatletter
\usepackage{breqn}
\usepackage{pdflscape}
\usepackage{afterpage}
\usepackage{capt-of}
\usepackage{dsfont}
\usepackage{gensymb}
\usepackage[mathscr]{euscript} 
\newcolumntype{x}[1]{>{\raggedright\arraybackslash\hspace{0pt}}p{#1}}
\newtheorem{problem}{Problem}
\newtheorem{assumption}{Assumption}

\newcommand{\iter}[1]{\ensuremath{\langle #1 \rangle}}
\newcommand{\itrpt}[3]{\ensuremath{\bm{#1}^{\iter{#2}}_{#3}}}
\newcommand{\itrval}[3]{\ensuremath{#1^{\iter{#2}}_{#3}}}

\newcommand{\idxpt}[3]{\ensuremath{\bm{#1}^{(#2)}_{#3}}}
\newcommand{\idxptT}[3]{\ensuremath{\bm{#1}^{(#2)\top}_{#3}}}
\newcommand{\idxval}[3]{\ensuremath{#1^{(#2)}_{#3}}}

\begin{document}
\begin{frontmatter}

\title{Rethinking Physics-Informed Regression Beyond Training Loops and Bespoke Architectures\thanksref{footnoteinfo}} 

\thanks[footnoteinfo]{This work has been submitted to IFAC for possible publication and was supported by the UKRI EPSRC Grant EP/Z533816/1 ``Concurrent Learning and Control of Uncertain Large-Scale Phenomena''.}

\author[First]{Lorenzo Sabug, Jr.} 
\author[First,Second]{Eric Kerrigan}

\address[First]{Department of Electrical \& Electronic Engineering, Imperial College, SW7 2AZ London, UK (e-mail: l.sabug21@imperial.ac.uk).}
\address[Second]{Department of Aeronautics, Imperial College, \\SW7 2AZ London, UK (e-mail: e.kerrigan@imperial.ac.uk)}

\begin{abstract}
. We revisit the problem of physics-informed regression, and propose a method that directly computes the state at the prediction point, simultaneously with the derivative and curvature information of the existing samples. We frame each prediction as a constrained optimisation problem, leveraging multivariate Taylor series expansions and explicitly enforcing physical laws. Each individual query can be processed with low computational cost without any pre- or re-training, in contrast to global function approximator-based solutions such as neural networks. Our comparative benchmarks on a reaction-diffusion system show competitive predictive accuracy relative to a neural network-based solution, while completely eliminating the need for long training loops, and remaining robust to changes in the sampling layout.
\end{abstract}

\begin{keyword}
Physics-informed machine learning, partial differential equations, regression, data-driven methods, nonlinear optimisation
\end{keyword}

\end{frontmatter}


\section{Introduction}

In many areas of science and engineering, such as those involving fluid dynamics \citep{Sharma2023-ta}, propagation phenomena \citep{Deng2025-nw}, and material diffusion/transport \citep{Xue2024-vo}, we encounter spatio-temporal processes \citep{Mercieca2016-ml}, whose state evolution is governed by partial differential equations (PDEs). While the state evolution is continuous in space and time, the resulting field cannot be fully observed due to physical and/or practical constraints. Instead, we only sense them at discrete spatial locations, and the measurements can be used to reconstruct (and eventually, also forecast) the states in the entire domain of interest. Because sensor data are only available at discrete locations and are often sparse or incomplete, reconstructing the full field is an ill-posed, with infinitely many fields compatible with the data. Therefore, regularisation using \textit{a priori} physical knowledge (e.g., using PDE information) can constrain the possible solutions and improve the fidelity of the state reconstruction with respect to the underlying process. This philosophy of incorporating physical laws into data-driven modelling has been the driving force behind a class of powerful approaches known by the broad category of ``physics-informed machine learning'' (PIML) methods \citep{Karniadakis2021-bd,Meng2025-si}.

Under the general coverage of PIML, we consider in this paper the so-called physics-informed regression models. Among the most prominent approaches in the literature are physics-informed neural network (PINN)-based approaches and physics-informed Gaussian processes (PIGPs) \citep{Cross2024-ic,Hanuka2021-ta}. Despite their wide popularity, both approaches suffer from the following drawbacks:
\begin{itemize}
    \item PINNs are trained to jointly minimise data and PDE-related losses in a single, non-convex optimisation problem, leading to unstable training and requiring many epochs to converge \citep{Wang2021-uy}. Furthermore, the performance of the PINN (or neural networks in general) is strongly influenced by its hyperparameters, architecture, and their interactions with the training data layout \citep{Cuomo2022-yz}. As a result, every new data layout and requirement might necessitate a redesign of the architecture and hyperparameters of the PINN. This makes the PINN setup process more of a \textit{bespoke engineering task} (requiring trial-and-error, each involving full retraining) rather than an out-of-the-box regression method.
    \item PIGPs alleviate some engineering setup problems and are easier to deploy than PINNs; however, they are not without their own issues. For instance, only linear PDEs can be analytically embedded within the GP prior \citep{Pfortner2022-di,Chen2021-cu}, and their resulting fidelity to the underlying process also depends on the choice of kernel. Lastly, their cubic computational complexity limits their scalability to larger datasets. While there have already been computational workarounds built on local, sparse, or inducing-point approximations, these trade off accuracy for speed and make the PDE embedding even more difficult.
    \item \textit{Both} of these method classes suffer from the implicitness or opaqueness of derivative information \citep{De_Ryck2024-ce}, which is a drawback for interpretability. The derivatives in PINNs are buried in the neuron weights and offsets, while for PIGPs, they are encoded in the kernel hyperparameters.
\end{itemize}

In this paper, we revisit physics-informed regression by proposing a much simpler, completely interpretable, and optimisation-based formulation. Instead of training a parametric surrogate function that implicitly satisfies \textit{a priori} physical knowledge, we frame each prediction as a constrained optimisation problem that is solved independently for each prediction point, explicitly enforcing consistency with the data and the physics. The data consistency is derived from Taylor series relations between the prediction point and the existing samples, which are encoded as linear inequality constraints. On the other hand, the agreement with physical laws is handled by using PDE-specific equality constraints. Our method, dubbed Direct Constraints-Based Regression (DCBR), offers the following advantages:

\begin{itemize}
    \item Being fully non-parametric, we do not assume any library of pre-defined basis functions on the data, skipping the need for residual minimisation-based fitting.
    \item Instead of hiding the derivatives in model parametrisations, DCBR explicitly computes the field variables alongside the state prediction as part of the decision variables in an optimisation problem, enabling direct interpretability in predictive regression.
    \item DCBR predicts the state value directly from the data set using constrained optimisation, eliminating the need for model pre-training, retraining and architecture tuning.
    \item The enforcement of \textit{a priori} known PDEs is explicit and transparent: they are simply enforced as equality constraints on the state values and derivatives.
\end{itemize}

\noindent Our comparative benchmarks on a reaction-diffusion system show that our method achieves competitive accuracy, and sometimes performs better than a PINN-based solution, without the need for long pre- or re-training loops.

This paper is organised into six sections. Section~\ref{sec:preliminaries} introduces the conventions and nomenclature, followed by the introduction of our proposed technique in Section~\ref{sec:dcbr}. A comparative case study in a reaction-diffusion system is discussed in Section~\ref{sec:case-study}, with the results analysed in Section~\ref{sec:results-discussion}. Finally, we lay out our concluding remarks in Section~\ref{sec:conclusions}.

\section{Preliminaries}
\label{sec:preliminaries}

We denote vectors as bold lowercase (e.g., $\bm{x}$, $\bm{\xi}$), matrices as bold uppercase (e.g., $\bm{H}$), and $\|\cdot\|$ as the Euclidean norm when operated on a vector. Furthermore, we mark a member $i$ of a set by the superscript $(\cdot)^{(i)}$, and we denote a quantity/set that changes with the iteration $k$ through the superscript $(\cdot)^{\langle k \rangle}$. We denote any point in 2D space as $\bm{p} \doteq [p_1 ~ p_2]^\top \in \mathcal{R}$ where $\mathcal{R} \subset \mathbb{R}^2$ is a closed and compact two-dimensional\footnote{While we consider 2D space in this paper for simplicity of discussion, the proposed method can be easily extended to higher dimensions.} region, and $t \in \mathbb{R}_+$ as the continuous time. Now we consider a spatio-temporal process on the domain $\mathcal{X} \doteq \mathcal{R} \times \mathbb{R}_+$, which is described by a partial differential equation (PDE)
\begin{equation}
    \mathcal{F}(u(\bm{x}), \mathcal{D} u(\bm{x}),\mathcal{D}^2u(\bm{x}),\ldots) = 0
\label{eqn:general-pde}
\end{equation}

\noindent where $\bm{x} \doteq [\bm{p}^\top ~ t]^\top$ is the spatio-temporal variable, and $\mathcal{D}$ is the differential operator acting on all components of the state $u$. Suppose now that we have collected a data set $\itrpt{X}{n}{}$ composed of $n$ samples
\[
    \itrpt{X}{n}{} \doteq \left\{(\idxpt{x}{i}{}, \idxval{u}{i}{})\right\}_{i=1}^n 
\]
\noindent where $\idxpt{x}{i}{} = [\idxpt{p}{i}{}{}^\top ~ \idxval{t}{i}{}]^\top \in \mathcal{X}$ is the sampled point, and $\idxval{u}{i}{}$ is the corresponding state. First, we state an assumption regarding the PDE

\begin{assumption}[Smooth PDE solution]
    There exists a solution $u$ to the PDE \eqref{eqn:general-pde} such that   
    \[
        u \in C^2(\mathcal{X}).
    \]
    Furthermore, along a segment from any $\idxpt{x}{i}{}$ to any $\bm{x}' \in \mathcal{X}$, the third and fourth directional derivatives of $u$ are bounded.
\label{assm:smooth-pde}
\end{assumption}

Our second assumption pertains to our sampling, stated as follows:
\begin{assumption}[Noise-free measurements]
    For every sample $(\idxpt{x}{i}{}, \idxval{u}{i}{}) \in \itrpt{X}{n}{}$,
    \[
        \idxval{u}{i}{} = u(\idxpt{x}{i}{}).
    \]
\label{assm:noise-free}
\end{assumption}
In laboratory and other controlled setups, Assumption~\ref{assm:noise-free} can be taken, and is convenient for discussion of our method; the treatment of finite additive noise is a subject of further work. Now we are faced with the problem of predicting the state~$u'$ at an unsampled spatio-temporal point $\bm{x}' = [\bm{p}'{^\top} ~ t']^\top$, which we refer to as the \textit{query point}. The problem is stated as follows:

\begin{problem}[Physics-consistent pointwise prediction]
    Given Assumptions~\ref{assm:smooth-pde} and \ref{assm:noise-free}, as well as the training data set $\itrpt{X}{n}{}$, governing PDE $\mathcal{F}$, and a query point $\bm{x}' \in \mathcal{X}$, determine the state value $u' \doteq u(\bm{x}')$ such that
    \[
        u' \in \mathcal{S}(\bm{x}'; \mathcal{F}, \itrpt{X}{n}{})
    \]
    \noindent where $\mathcal{S}$ is the set of values at $\bm{x}$ consistent with $\mathcal{F}$ and~$\itrpt{X}{n}{}$.
\label{prob:physics-prediction}
\end{problem}

For the purposes of this paper, we refer to estimation as \textit{interpolation} when $\min_i \idxval{t}{i}{} \leq t' \leq \max_i \idxval{t}{i}{}$, and as \textit{forecast} otherwise. Furthermore, we use $\bm{\xi}_i \doteq \bm{x}' - \idxpt{x}{i}{}$, and $\bm{\xi}_{ij} \doteq \idxpt{x}{j}{} - \idxpt{x}{i}{}$. We denote $\partial_1 u, \partial_2 u, \partial_t u$ as the derivatives of $u$ w.r.t.\ $p_1, p_2, t$, respectively, and $\partial_1^2 u, \partial_2^2 u, \partial_t^2 u$ the second derivatives. Lastly, $\bm{g}(\bm{x}) \doteq \mathcal{D} u(\bm{x})$ is the spatio-temporal gradient, $\bm{H}(\bm{x}) \doteq \mathcal{D}^2u(\bm{x})$ the Hessian, and we use the shorthand notations $\idxpt{g}{i}{} = \bm{g}(\idxpt{x}{i}{})$ and $\idxpt{H}{i}{} = \bm{H}(\idxpt{x}{i}{})$.

\section{Direct Constraints-Based Regression~(DCBR)}
\label{sec:dcbr}

We address Problem~\ref{prob:physics-prediction} by framing the prediction problem as an optimisation and feasibility problem, enforcing Taylor-based consistency and physical laws as constraints. 

\subsection{Query-to-sample Taylor constraints}
Let us consider the multivariate Taylor series expansion with integral remainder, considering an existing sample $(\idxpt{x}{i}{},\idxval{u}{i}{}) \in \itrpt{X}{n}{}$ and the state value $u'$ at the prediction point $\bm{x}'$: 
\begin{equation}
    u' \approx \idxval{u}{i}{} + \idxptT{g}{i}{} \bm{\xi}_i + \int_0^1 (1-s)\bm{\xi}^\top_i\bm{H}(\idxpt{x}{i}{} + s\bm{\xi}_i) \bm{\xi}_i ds
    \label{eqn:taylor-series}
\end{equation}

A particular challenge with using the Taylor series expansion is that we are only provided with a zeroth-order data set (i.e., only the sampled locations/times and corresponding state values) and we lack gradient and curvature information. We tackle this by (1) first exploiting a quadrature rule to approximate the integral, (2) embedding the gradient kinematic relations, and (3) writing similar Taylor series relations among different samples to enforce consistency, and thereby constrain the gradient and curvature estimates.

Numerically approximating the integral at \eqref{eqn:taylor-series} using a 3-point Gauss-Lobatto quadrature, we have
{\small
\[
    u' \approx \idxval{u}{i}{} + \idxptT{g}{i}{} \bm{\xi}_i + \sum_{k=0}^2w_k (1 - s_k) \bm{\xi}^\top_i \bm{H}(\idxpt{x}{i}{} +  s_k\bm{\xi}_i)\bm{\xi}_i.
\]}

With $s_k = \{0,0.5,1\}$ and the corresponding weights $w_k = \{\frac{1}{6},\frac{4}{6},\frac{1}{6}\}$ for $k=\{0,1,2\}$, respectively, we arrive at
{\small
\begin{equation}
    u' \approx \idxval{u}{i}{} + \idxptT{g}{i}{} \bm{\xi}_i + \frac{1}{6}\bm{\xi}^\top_i\idxpt{H}{i}{}\bm{\xi}_i + \frac{1}{3}\bm{\xi}^\top_i\bm{H}(\idxpt{x}{i}{}+\frac{1}{2}\bm{\xi}_i)\bm{\xi}_i,
\label{eqn:taylor-series-quadrature}
\end{equation}}
\noindent where $\bm{H}(\idxpt{x}{i}{} + \frac{1}{2}\bm{\xi}_i)$ we now denote as $\bm{H}_m$, which is the Hessian at the midpoint from $\idxpt{x}{i}{}$ to $\bm{x}'$.

Invoking the kinematic relation between the corresponding gradients, and using a 3-point Gauss-Lobatto quadrature again, we get
{\small
\begin{align}
    \bm{g}' - \idxpt{g}{i}{} &= \int_0^1 (\idxpt{H}{i}{} + s\bm{\xi}_i) \bm{\xi}_i ds \\
    ~ &\approx \tfrac{1}{6} \idxpt{H}{i}{}\bm{\xi}_i + \tfrac{2}{3}\bm{H}_m \bm{\xi}_i + \tfrac{1}{6}\bm{H}' \bm{\xi}_i,
\end{align}}
\noindent and now introducing the approximation error,
{\small
\[
    \bm{H}_m \bm{\xi}_i = \tfrac{3}{2}(\bm{g}' - \idxpt{g}{i}{}) - \tfrac{1}{4}(\bm{H}' + \idxpt{H}{i}{})\bm{\xi}_i + \epsilon_g
\]}

\noindent where $|\epsilon_g| \leq (\epsilon_g' + \epsilon^{(i)}_g) r_i^4$ is the quadrature remainder for the gradient, with $r_i \doteq \|\bm{\xi}_i\|$. Plugging into \eqref{eqn:taylor-series-quadrature} and also putting the corresponding error term $\epsilon_Q$,
{\small
\begin{multline}
    u' = \idxval{u}{i}{} + \idxptT{g}{i}{} \bm{\xi}_i + \frac{1}{6}\bm{\xi}^\top_i\idxpt{H}{i}{}\bm{\xi}_i ~ + ~\\ \frac{1}{3} \bm{\xi}^\top_i \left( \frac{3}{2}(\bm{g}' - \idxpt{g}{i}{}) - \frac{1}{4}(\bm{H}' + \idxpt{H}{i}{})\bm{\xi}_i + \epsilon_g\right) + \epsilon_Q
\end{multline}}

\noindent where $|\epsilon_Q| \leq \idxval{\epsilon}{i}{Q} r_i^4$ is the quadrature remainder for $u'$. Now we denote $\bm{h} = \text{vec}(\bm{H})$, $\bm{\lambda}_i = \text{vec}(\bm{\xi}_i\bm{\xi}_i^\top)$, $\bm{\kappa}' \doteq [{\bm{g}'}^\top ~ {\bm{h}'}^\top ~ \epsilon'_g ~ \epsilon'_Q]^\top$, $\idxpt{\kappa}{i}{} \doteq [\idxptT{g}{i}{} ~ \idxptT{h}{i}{} ~ \idxval{\epsilon}{i}{g} ~ \idxval{\epsilon}{i}{Q}]^\top$, and absorb the constant $\frac{1}{3}$ into $\epsilon'_g, \idxval{\epsilon}{i}{g}, \idxval{\epsilon}{i}{Q}$. Collecting and reorganising to matrix form, we now have the inequalities


{\scriptsize
\begin{equation}
    \begin{bmatrix}
    1 & -\frac{1}{2}\bm{\xi}^\top_i & \frac{1}{12}\bm{\lambda}_i & -r_i^5 & 0 & -\frac{1}{2}\bm{\xi}^\top_i & -\frac{1}{12}\bm{\lambda}_i & -r_i^5 & -r_i^4 \\
    -1 & \frac{1}{2}\bm{\xi}^\top_i & -\frac{1}{12}\bm{\lambda}_i & -r_i^5 & 0 & \frac{1}{2}\bm{\xi}^\top_i & \frac{1}{12}\bm{\lambda}_i & -r_i^5 & -r_i^4
    \end{bmatrix}
    \begin{bmatrix}
        u' \\ \bm{\kappa}' \\ \idxpt{\kappa}{i}{}
    \end{bmatrix}
    \leq
    \begin{bmatrix}
        \idxval{u}{i}{} \\ -\idxval{u}{i}{}
    \end{bmatrix}
\end{equation}}

With $\bm{\xi}'_i \doteq \idxpt{x}{i}{}-\bm{x}'$, $r'_i \doteq \|\bm{\xi}'_i\|$, and $\bm{\lambda}'_i \doteq \text{vec}(\bm{\xi}'_i\bm{\xi}'_i{}^\top)$, we now add the backwards Taylor series expansion (expressing $\idxval{u}{i}{}$ w.r.t.  $u', \bm{g}', \bm{h}'$ in the same form as \eqref{eqn:taylor-series}), such that we get
{\scriptsize
\begin{multline}
    \begin{bmatrix}
    1 & -\frac{1}{2}\bm{\xi}^\top_i & \frac{1}{12}\bm{\lambda}_i & -r_i^5 & 0 & -\frac{1}{2}\bm{\xi}^\top_i & -\frac{1}{12}\bm{\lambda}_i & -r_i^5 & -r_i^4 \\
    -1 & \frac{1}{2}\bm{\xi}^\top_i & -\frac{1}{12}\bm{\lambda}_i & -r_i^5 & 0 & \frac{1}{2}\bm{\xi}^\top_i & \frac{1}{12}\bm{\lambda}_i & -r_i^5 & -r_i^4 \\    
    -1 & -\frac{1}{2}{\bm{\xi}'_i}^\top & -\frac{1}{12}\bm{\lambda}'_i & -{r'_i}^5 & -{r'_i}^4 & -\frac{1}{2}{\bm{\xi}'_i}^\top & \frac{1}{12}\bm{\lambda}'_i & -{r'_i}^5 & 0 \\
    1 & \frac{1}{2}{\bm{\xi}'_i}^\top & \frac{1}{12}\bm{\lambda}'_i & -{r'_i}^5 & -{r'_i}^4 & \frac{1}{2}{\bm{\xi}'_i}^\top & -\frac{1}{12}\bm{\lambda}'_i & -{r'_i}^5 & 0
    \end{bmatrix}
    \begin{bmatrix}
        u' \\ \bm{\kappa}' \\ \idxpt{\kappa}{i}{}
    \end{bmatrix}\\
    \leq
    \begin{bmatrix}
        \idxval{u}{i}{} \\ -\idxval{u}{i}{} \\ -\idxval{u}{i}{} \\ \idxval{u}{i}{}
    \end{bmatrix}
    \label{eqn:stencil-prediction}
\end{multline}}

For the sake of brevity, we denote the first column as $\bm{M}$, the next four columns as $\bm{L}_{i}$, and the right four as $\bm{R}_{i}$. Denoting furthermore that $\bm{\delta}_i = \begin{bmatrix}\idxval{u}{i}{} & -\idxval{u}{i}{} & -\idxval{u}{i}{} & \idxval{u}{i}{}\end{bmatrix}^\top$, we now have
\[
    \begin{bmatrix}
        \bm{M} ~ \bm{L}_i ~ \bm{R}_i
    \end{bmatrix}
    \begin{bmatrix}
        u' \\ \bm{\kappa}' \\ \idxpt{\kappa}{i}{}
    \end{bmatrix}
    \leq
    \bm{\delta}_i.
\]

\subsection{Sample-to-sample Taylor constraints}
We can also impose Taylor series-based cross-consistency constraints between two existing samples $(\idxpt{x}{i}{},\idxval{u}{i}{})$, $(\idxpt{x}{j}{},\idxval{u}{j}{})$. Denoting $\bm{\lambda}_{ij} = \text{vec}(\bm{\xi}_{ij}\bm{\xi}_{ij}^\top)$, this leads us to
{\small
\begin{multline}
    \begin{bmatrix}        
    -\frac{1}{2}\bm{\xi}^\top_{ji} & \frac{1}{12}\bm{\lambda}_{ji} & -r_{ji}^5 & 0 & -\frac{1}{2}\bm{\xi}^\top_{ji} & -\frac{1}{12}\bm{\lambda}_{ji} & -r_{ji}^5 & -r_{ji}^4 \\
    \frac{1}{2}\bm{\xi}^\top_{ji} & -\frac{1}{12}\bm{\lambda}_{ji} & -r_{ji}^5 & 0 & \frac{1}{2}\bm{\xi}^\top_{ji} & \frac{1}{12}\bm{\lambda}_{ji} & -r_{ji}^5 & -r_{ji}^4 \\
    -\frac{1}{2}\bm{\xi}^\top_{ij} & -\frac{1}{12}\bm{\lambda}_{ij} & -r_{ij}^5 & -r_{ij}^4 & -\frac{1}{2}\bm{\xi}^\top_{ij} & \frac{1}{12}\bm{\lambda}_{ij} & -r_{ij}^5 & 0 \\
    \frac{1}{2}\bm{\xi}^\top_{ij} & \frac{1}{12}\bm{\lambda}_{ij} & -r_{ij}^5 & -r_{ij}^4 & \frac{1}{2}\bm{\xi}^\top_{ij} & -\frac{1}{12}\bm{\lambda}_{ij} & -r_{ij}^5 & 0
    \end{bmatrix}
    \begin{bmatrix}
        \idxpt{\kappa}{i}{} \\ \idxpt{\kappa}{j}{}
    \end{bmatrix} \\
    \leq
    \begin{bmatrix}
        \idxval{u}{j}{} - \idxval{u}{i}{} \\ \idxval{u}{i}{} - \idxval{u}{j}{} \\ \idxval{u}{i}{} - \idxval{u}{j}{} \\ \idxval{u}{j}{} - \idxval{u}{i}{} 
    \end{bmatrix},
    \label{eqn:stencil-2samples}
\end{multline}}

\noindent Denoting the matrix to the right side of the inequality as~$\bm{\delta}_{ij}$, we can simply shorten this to
\[
    \begin{bmatrix}
        \bm{L}_{ij} & \bm{R}_{ij}
    \end{bmatrix}
    \begin{bmatrix}
        \idxpt{\kappa}{i}{} \\ \idxpt{\kappa}{j}{}
    \end{bmatrix}
    \leq
    \bm{\delta}_{ij}.
\]

\subsection{Aggregated constraints and optimisation problem}

Putting all the Taylor series relations together, we have
{\small
\begin{equation}
    \underbrace{
    \left[
    \begin{array}{c}
    \begin{matrix}
        \bm{M} & \bm{L}_1 & \bm{R}_1 & \bm{0} & \cdots & \bm{0} \\
        \bm{M} & \bm{L}_2 & \bm{0}   & \bm{R}_2 & ~ & \bm{0} \\
        \vdots & \vdots   & \vdots   & ~   & \ddots & \vdots \\
        \bm{M} & \bm{L}_n & \bm{0}   & \bm{0} & \cdots & \bm{R}_n \\
    \end{matrix} \\
    \midrule
    \begin{matrix}
        \bm{0} & \bm{0} & \bm{L}_{12} & \bm{R}_{12} & \cdots & \bm{0} \\
        \vdots & \vdots & \vdots      & \vdots      & \ddots & \vdots \\
        \bm{0} & \bm{0} & \bm{L}_{1n} & \bm{0}      & \cdots & \bm{R}_{1n} \\
    \end{matrix} \\
    \midrule
    \begin{matrix}
        \vdots
    \end{matrix} \\ 
    \midrule
    \begin{matrix}
        \bm{0} & \bm{0} & \cdots & \bm{0}      & \bm{L}_{n-1,n} & \bm{R}_{n-1,n} \\
    \end{matrix}
    \end{array}
    \right]
    }_{\bm{A}}
    \underbrace{
    \begin{bmatrix}
        u' \\ \bm{\kappa}' \\ \idxpt{\kappa}{1}{} \\ \idxpt{\kappa}{2}{} \\ \idxpt{\kappa}{3}{} \\ \vdots \\ \idxpt{\kappa}{n}{}
    \end{bmatrix}
    }_{\bm{\theta}}
    \leq
    \underbrace{
    \left[
    \begin{array}{c}
    \begin{matrix}
        \bm{\delta}_1 \\ \bm{\delta}_2 \\ \vdots \\ \bm{\delta}_n 
    \end{matrix} \\
    \midrule
    \begin{matrix}
        \bm{\delta}_{12} \\ \vdots \\ \bm{\delta}_{1n} 
    \end{matrix} \\
    \midrule
    \vdots \\
    \midrule
    \bm{\delta}_{n-1,n}\\
    \end{array}
    \right]
    }_{\bm{b}}
\end{equation}}

For brevity, we denote $\bm{\omega} \doteq [\bm{g}^\top ~ \bm{h}^\top]^\top$ as the field variables such that $\bm{\kappa} = [\bm{\omega}^\top ~ \epsilon_g ~ \epsilon_Q]^\top$. To set up our optimisation problem, we now minimise the L2 penalty involving our slack variables, as follows:
\begin{subequations}
\begin{alignat}{4}
    \min_{\bm{\theta}} ~& {\epsilon'_g}^2 + {\epsilon'_Q}^2 + \sum_i \epsilon_g^{(i)2}  + \sum_i \epsilon_Q^{(i)2} + \rho \bm{\theta}^\top \bm{\theta} \label{eqn:nlp-objective} \\
    \text{s.t.} &~\bm{A}\bm{\theta} \leq \bm{b} \label{eqn:nlp-inequalities} \\
    ~& \epsilon'_g, \epsilon'_Q, \idxval{\epsilon}{i}{g}, \idxval{\epsilon}{i}{Q} \geq 0 \\
    ~& \bm{\omega}' \in \mathcal{P} \label{eqn:nlp-equalities-1}\\
    ~& \idxpt{\omega}{i}{} \in \mathcal{P}, \forall i \in \{1, \ldots, n\}. \label{eqn:nlp-equalities-2}
\end{alignat}
\end{subequations}
\noindent where $\rho$ is a small value for regularisation (in this paper, $\rho=1\times10^{-9}$). Our decision variable $\bm{\theta}$ contains the predicted state, as well as the derivative estimates at the query point $\bm{x}'$ and the existing samples $\idxpt{x}{i}{}$. $\mathcal{P}$ is the set of $\bm{\omega}$ that satisfy the PDE. Therefore, \eqref{eqn:nlp-equalities-1}--\eqref{eqn:nlp-equalities-2} are implemented as equalities relating the (first, or higher-/mixed-) derivative terms of a point with each other, which, depending on the PDE considered, translate to linear or nonlinear equalities. In principle, \eqref{eqn:nlp-objective}--\eqref{eqn:nlp-equalities-2} finds the smallest error bounds consistent with the data and the physics (such the Taylor series inequalities \eqref{eqn:nlp-inequalities} and the PDE \eqref{eqn:nlp-equalities-1}--\eqref{eqn:nlp-equalities-2} are satisfied), applying Set Membership principles \citep{MILANESE2004957} to the physics-informed regression problem.





\section{Case Study}
\label{sec:case-study}

\begin{figure*}
    \centering
    \includegraphics[width=0.625\linewidth]{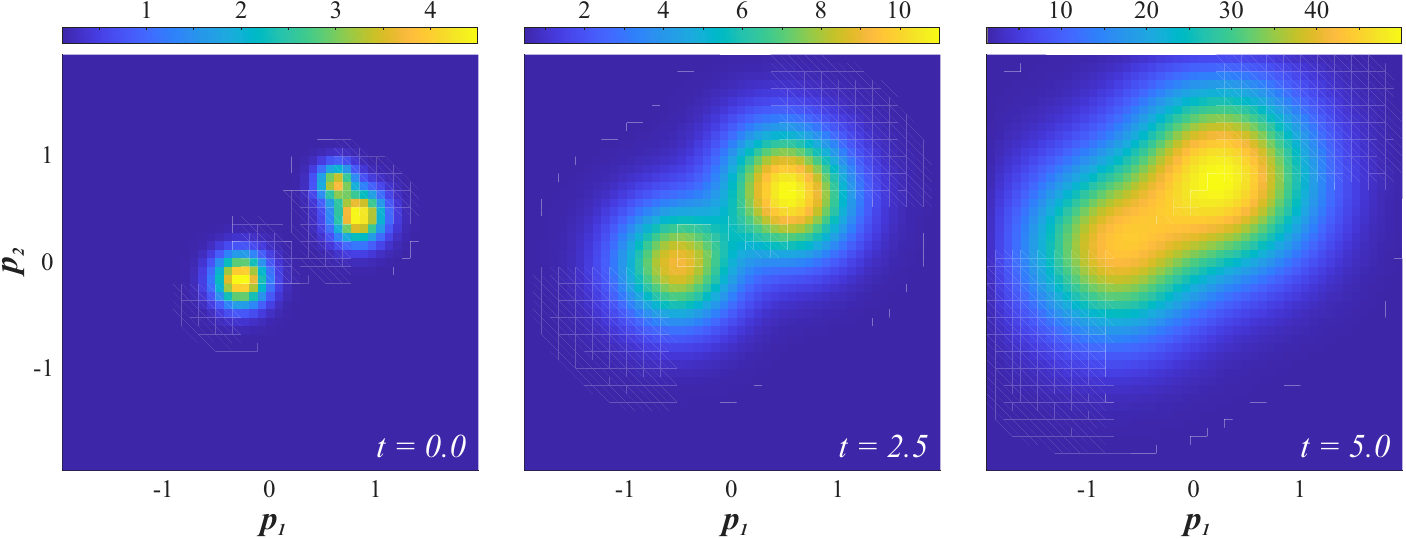}
    \caption{Chronological evolution snapshots of the reaction-diffusion system considered in this paper}
    \label{fig:rds-snapshots}
\end{figure*}

We consider a reaction-diffusion system (RDS), described by the PDE
\[
    \partial_t u = \nu\nabla^2u + \alpha u - \beta u^2 + \bm{w} \cdot \nabla u
\]
\noindent with coefficients $\nu=0.02$, $\alpha=1.0$, $\beta=0.008$, and $\bm{w} = \begin{bmatrix}0.1 & -0.06\end{bmatrix}^\top$. This PDE roughly approximates the dynamics of a forest fire with wind direction to the southwest. The initial condition (IC) is composed of three peaks and, solving the PDE in marching time, the peaks expand and eventually merge into one, as in Figure~\ref{fig:rds-snapshots}.

\subsection{Compared methods}
In this paper, we compare the following methods: 

\begin{itemize}
    \item Physics-informed neural network (PINN): the most widely adopted baseline for physics-informed machine learning. For this paper, we used the Python TensorFlow-based package \texttt{pinns-tf2}. We have used the following architecture: 3 input dimensions, 3 fully-connected hidden layers (128 neurons each layer, with ReLU activation), 1 linear output layer. The model was trained for 50000 epochs using the Adam optimizer (learning rate = 0.001), and with 10000 collocation points for the PDE loss minimisation.
    \item Direct constraints-based regression (DCBR, ours). We have set up the optimisation problem \eqref{eqn:nlp-objective}--\eqref{eqn:nlp-equalities-2}, with the PDE explicitly embedded as equality constraints. Recognising that $\partial_1 u'$, $\partial_2 u'$, $\partial_t u'$, $\partial_1^2u'$, and $\partial_2^2 u'$ are elements of $\bm{\omega}'$, we have the nonlinear equality \eqref{eqn:nlp-equalities-1} as
    \[
        \nu \partial^2_1 u' + \nu \partial^2_2u' + \alpha u' - \beta u'{}^2 + w_1 \partial_1 u' + w_2 \partial_2 u' - \partial_t u'= 0.
    \]
    \noindent However, for each existing sample, the state value $\idxval{u}{i}{}$ is already known; therefore, \eqref{eqn:nlp-equalities-2} conveniently reduces to a linear equality:
    \[
        \begin{bmatrix}
            -w_1 & -w_2 & 1 & \nu & \nu
        \end{bmatrix}
        \begin{bmatrix}
            \partial_1 \idxval{u}{i}{} \\
            \partial_2 \idxval{u}{i}{} \\
            \partial_t \idxval{u}{i}{} \\
            \partial_1^2 \idxval{u}{i}{} \\
            \partial_2^2 \idxval{u}{i}{}
        \end{bmatrix}
        = \alpha \idxval{u}{i}{} - \beta \idxval{u}{i}{} {}^2.
    \]
    Our proposed DCBR was implemented in Julia~1.12.1, and solved with \texttt{Ipopt} \citep{wachter2006ipopt}.
\end{itemize}

We recognise that there are other (more recent) methods that can embed physics into the regression process. One are the physics-informed Gaussian processes (PIGPs), which, however, are only natively able to handle linear PDEs \citep{Pfortner2022-di}. Therefore we did not include PIGPs for a direct comparison (whilst there are recent linearisation-based workarounds, we do not consider them in this paper). Another are operator networks \citep{Lu2021-ee}, however, we do not consider them because they will require training on data triples (each contains data samples on a \textit{fixed set} of locations, a query point outside of those locations, and a corresponding query value), which is completely different from PINN and DCBR data requirements.

\subsection{Data sets}

We have generated a high-resolution data set by running an RDS simulation and acquiring a 300$\times$300 sample grid per time step $k \in \{0, 1, \ldots, K\}$, which corresponds to the continuous time $t_k=0.1k$. The time step-related data set is considered the ``ground truth snapshot'' 
\[
    \itrpt{G}{k}{} = \left\{\left(\idxpt{p}{i_s}{}, t_k,u(\idxpt{p}{i_s}{},t_k)\right)\right\}
\]
\noindent where $\idxpt{p}{i_s}{}, i_s \in \{1,\ldots,300^2\}$ is a spatial point in the high-resolution grid, and $u(\idxpt{p}{i_s}{},t_k)$ is the corresponding PDE solution at $t=t_k$. From the ground truth snapshots, we build two training data sets:

\begin{itemize}
    \item Grid data set (\textbf{GRID}): for every $t_k$, we acquire a regular spatial grid of samples of different (low) resolutions (10$\times$10, 20$\times$20), resulting in a corresponding snapshot $\itrpt{S}{k}{G}$.
    \item Random data set (\textbf{RAND}): for every $t_k$, we get an independently and uniformly distributed set of samples, producing the snapshot $\itrpt{S}{k}{R}$. We have also tried different number of samples (100, 500) per time step.
\end{itemize}

\subsection{Tests and comparison metrics}

We use the compared methods PINN and DCBR to predict the (high-resolution) state snapshot at time $t=t_{k'}$,
\[
    \tilde{\bm{G}}^{\langle k' \rangle} = \left\{\left(\idxpt{p}{i_s}{}, t_{k'},\tilde{u}(\idxpt{p}{i_s}{},t_{k'})\right)\right\}
\]
and in this paper, we have two tests according to the provided data sets, and $t_{k'}$:

\begin{itemize}
    \item \textit{Spatial interpolation}: provided all (low resolution) snapshots $\itrpt{S}{0}{}, \ldots, \itrpt{S}{K}{}$, we compute $\tilde{\bm{G}}^{\langle k' \rangle}$, with $k' \in \{1,\ldots,K\}$,
    
    \item \textit{State forecast}: provided the (low resolution) snapshots $\itrpt{S}{0}{}, \ldots, \itrpt{S}{k}{}$, we compute the next snapshots $\tilde{\bm{G}}^{\langle k' \rangle}$, with $k' \in \{k+1,\ldots,K\}$.
\end{itemize}

We compared the accuracy of the methods using the (discretised) $L_2$ relative error, defined as
\begin{equation}
    \itrval{\epsilon}{k'}{} = \sqrt{\frac{\sum_{i_s} \left(\tilde{u}(\idxpt{p}{i_s}{},t_{k'}) - u(\idxpt{p}{i_s}{},t_{k'})\right)^2}{\sum_{i_s} u(\idxpt{p}{i_s}{},t_{k'})^2 + \varepsilon}}
\end{equation} \\
\noindent where $\varepsilon=1\times10^{-6}$ is a small constant for numerical stability. Furthermore, we compared their computational times, which we differentiate between:
\begin{itemize}
    \item \textit{Training}: time needed for the (PINN) model to be ready for the first query,
    \item \textit{Per-query}: average time to compute the prediction at a spatio-temporal point $\bm{x}'$.
\end{itemize}
All computations are done on an Ubuntu subsystem (Windows Subsystem for Linux) on a PC with Intel Core i9 (24~cores, up to 4.5~GHz), NVIDIA RTX 4500 Ada GPU, and 192 GB RAM. We have used the CUDA extensions for PINN training, while all DCBR-related computations are delegated in a single CPU thread.

\section{Results and discussion}
\label{sec:results-discussion}

\begin{figure*}
    \centering
    \includegraphics[width=\linewidth]{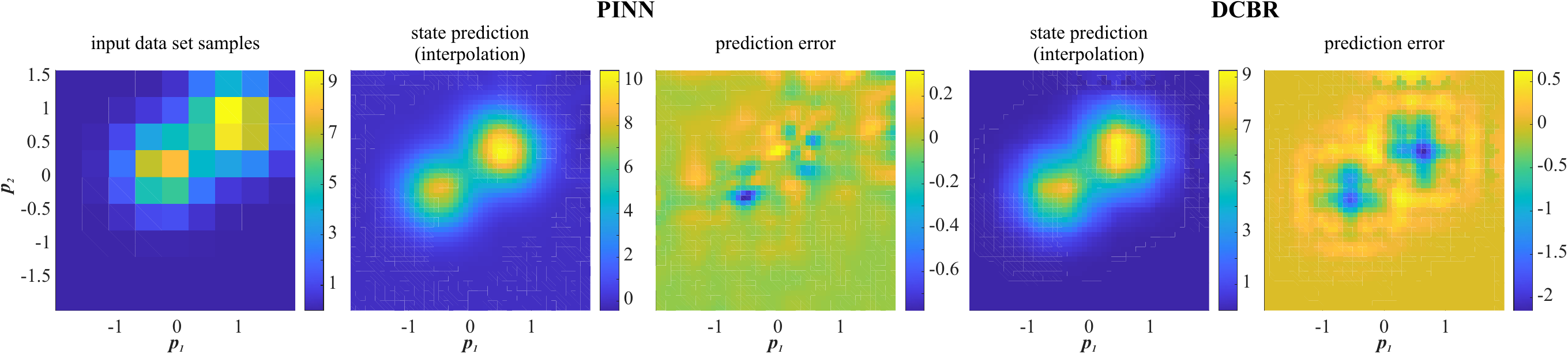}
    \caption{Case study (\textbf{GRID}, 10$\times$10): visual comparison and prediction errors at $t_k=2.5$\,s.}
    \label{fig:graphical-error-grid}
\end{figure*}
\begin{figure}
    \centering
    \includegraphics[width=\linewidth]{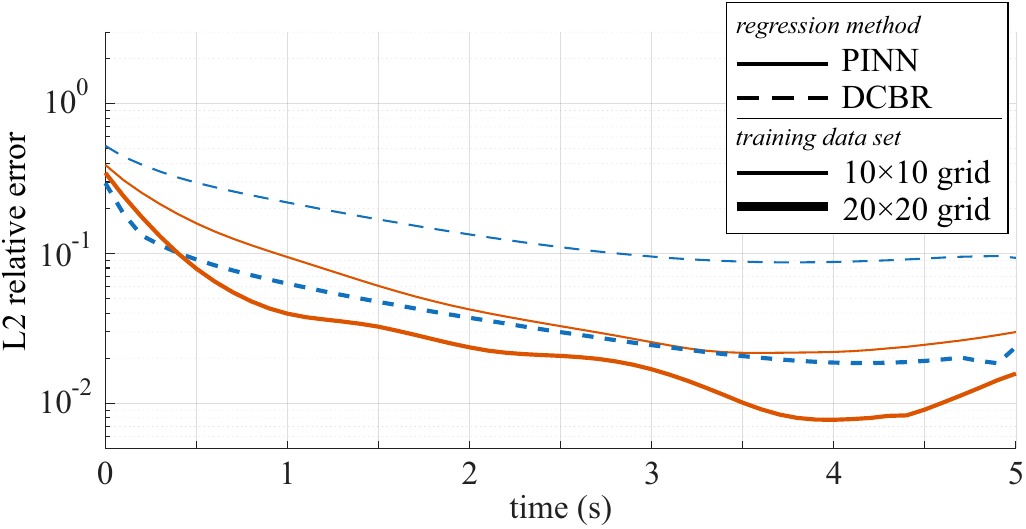}
    \caption{Case study (\textbf{GRID}): interpolation $L_2$ relative errors on different RDS snapshots}
    \label{fig:grid-l2-error-interpolation}
\end{figure}

In Figure~\ref{fig:graphical-error-grid}, we compare the resulting interpolated RDS state snapshots from a very low resolution (10$\times$10) grid data set at $t_k=2.5$\,s. In this ``upsampling'' exercise we observe that DCBR has a larger prediction error than PINNs, which is unsurprising due to two factors: (i) the input data set is composed of widely-spaced samples, and (ii) \eqref{eqn:nlp-objective}--\eqref{eqn:nlp-equalities-2} describes a local regressor using neighbouring data to predict the query point. This is in contrast with PINNs, where the training process shapes a globally coupled function representation, allowing data from distant locations to influence the predicted value at any given query point. This is confirmed by Figure~\ref{fig:grid-l2-error-interpolation}, in which PINN $L_2$ relative errors are much better than DCBR among all the snapshots considered. However, when we consider a (20$\times$20) finer grid data set, the $L_2$ relative errors have become very comparable across all snapshots.

\begin{figure}
    \centering
    \includegraphics[width=0.625\linewidth]{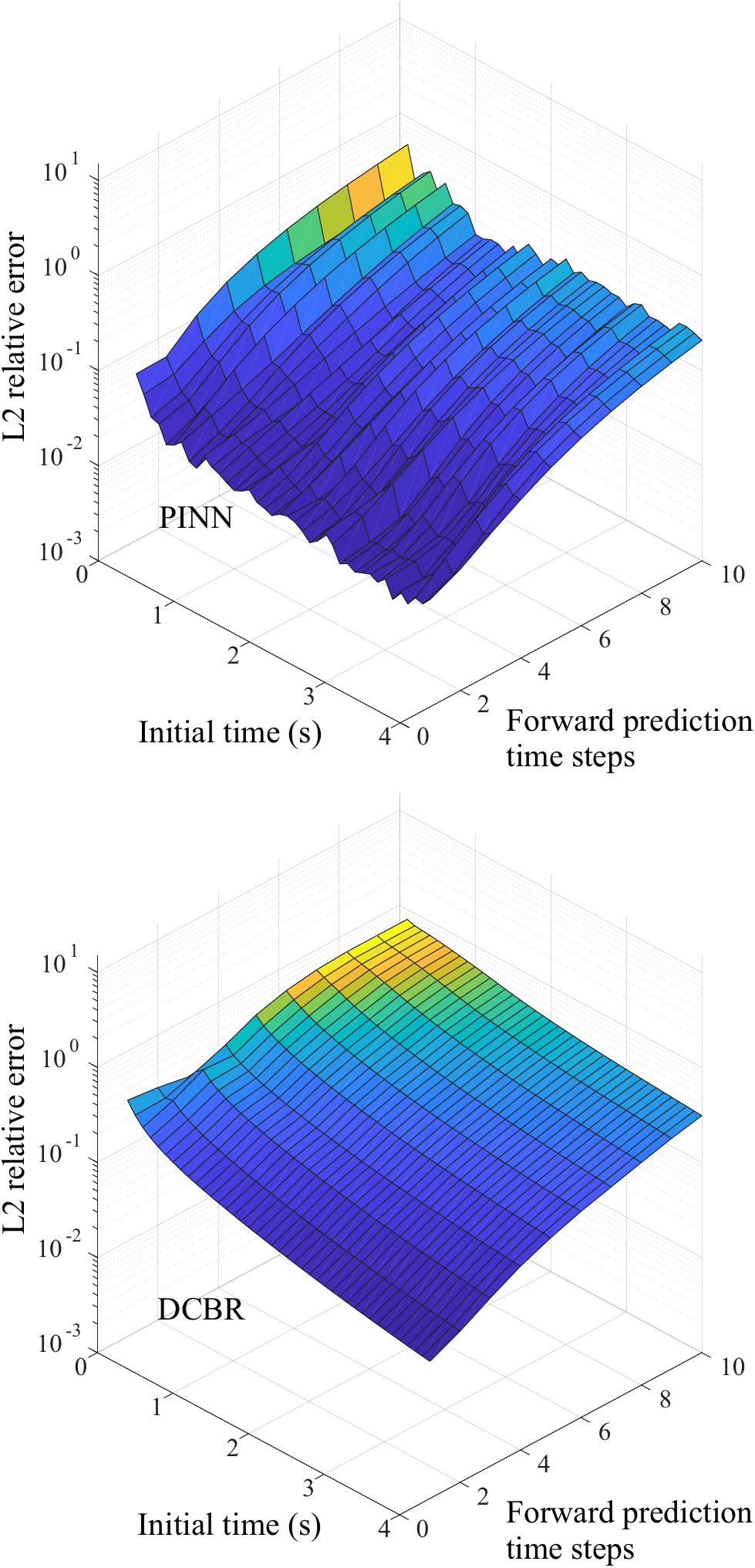}
    \caption{Case study (\textbf{GRID}): forecast $L_2$ relative errors with different ICs (20$\times$20 grid/timestep)}
    \label{fig:grid-l2-error-forecast}
\end{figure}

\begin{figure*}
    \centering
    \includegraphics[width=\linewidth]{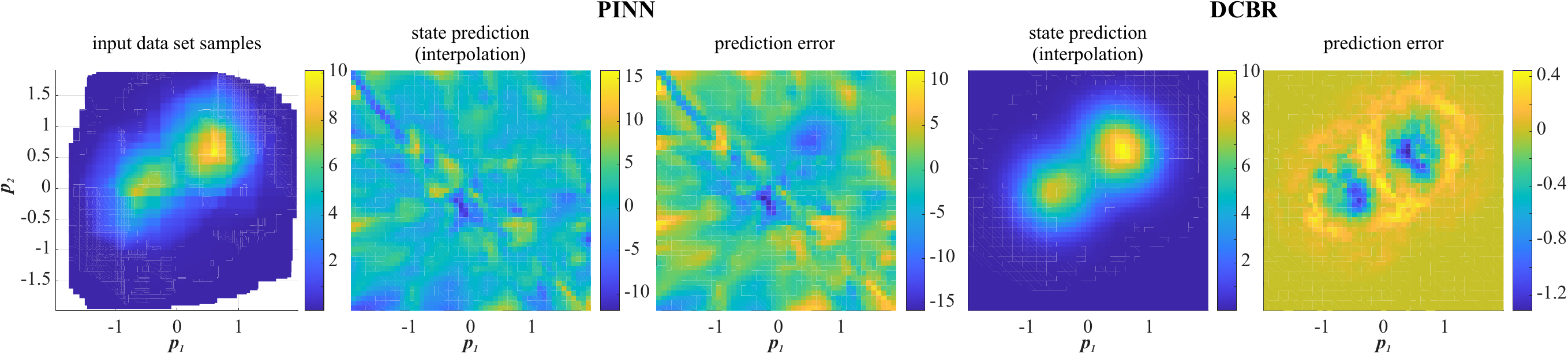}
    \caption{Case study (\textbf{RAND}, 100 samples): visual comparison and prediction errors at $t_k=2.5$\,s.}
    \label{fig:graphical-error-rand}
\end{figure*}

\begin{figure}
    \centering
    \includegraphics[width=\linewidth]{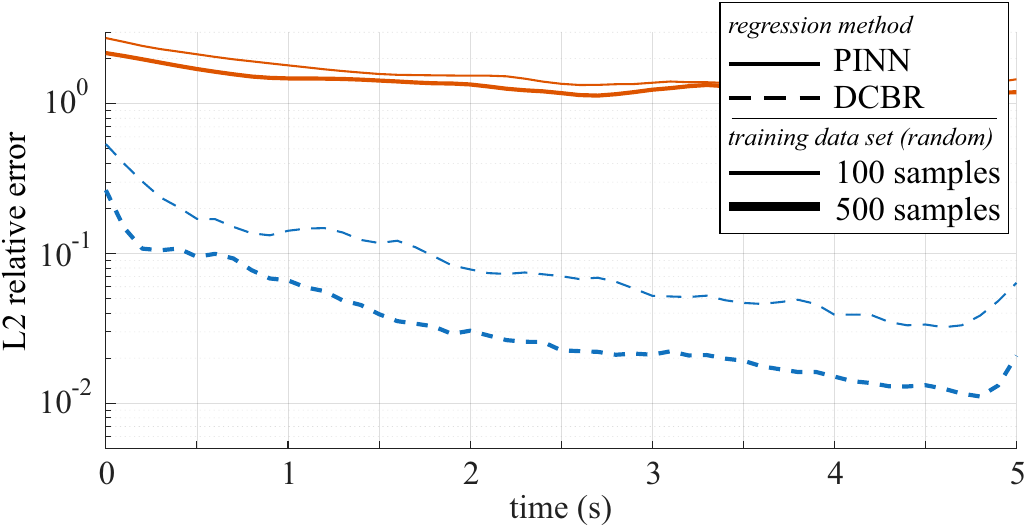}
    \caption{Case study (\textbf{RAND}): interpolation $L_2$ relative errors on different RDS snapshots}
    \label{fig:rand-l2-error-interpolation}
\end{figure}

\begin{figure}
    \centering
    \includegraphics[width=0.625\linewidth]{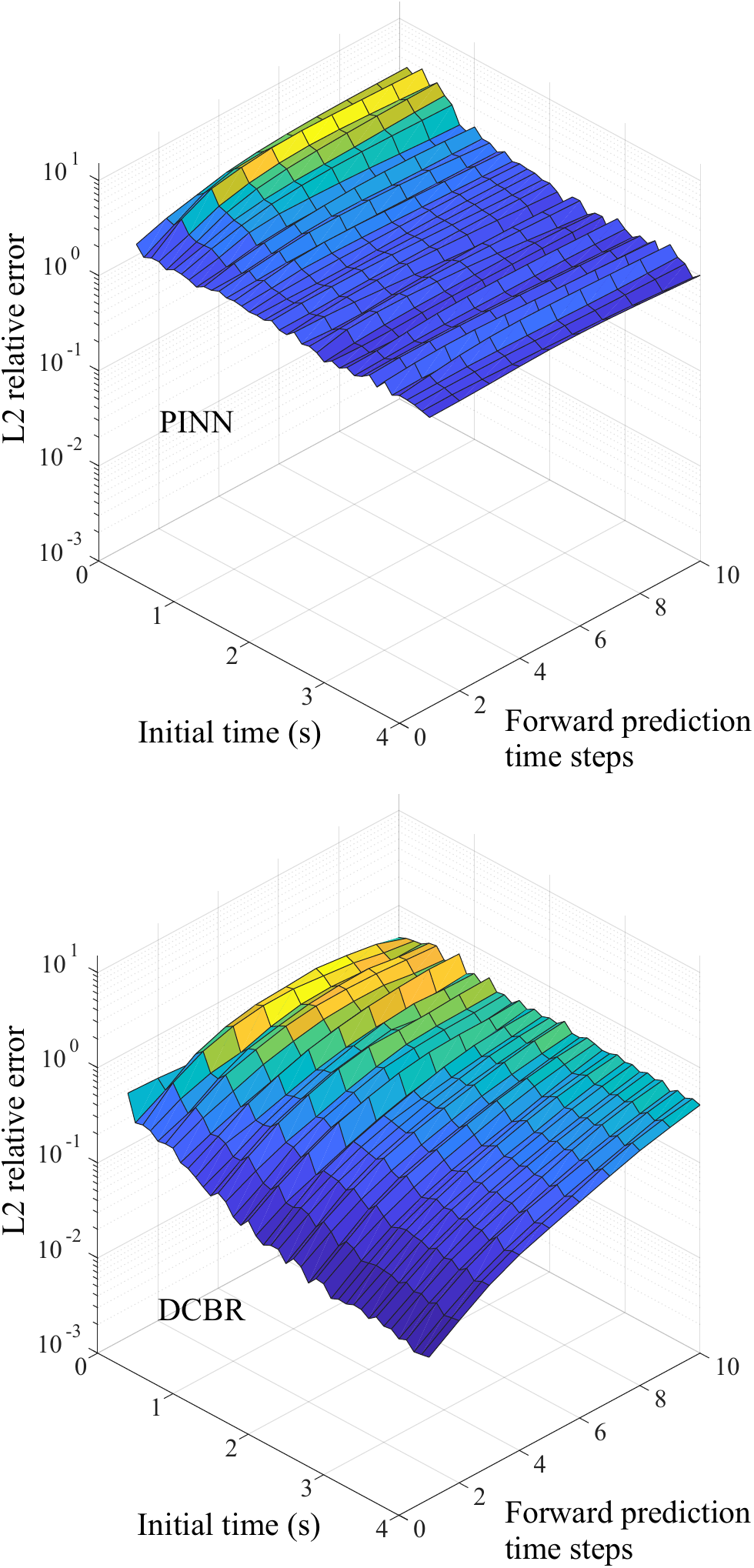}
    \caption{Case study (\textbf{RAND}): forecast $L_2$ relative errors with different ICs (500 samples/timestep)}
    \label{fig:rand-l2-error-forecast}
\end{figure}

In Figure~\ref{fig:grid-l2-error-forecast}, we show the forecast error metrics when starting from different ICs, given a 20$\times$20 grid data. When starting from an IC with higher-frequency content (e.g., IC is from the early stages of the event), the $L_2$ error evolution is worse than when given an IC with lower spatial frequency content. Furthermore, we note that the forecast error surfaces for PINN and DCBR are similar in magnitude, which is consistent with their similarity in accuracy for interpolation. However, we note that the PINN surface is much rougher because when we provide it with a different IC (and have to feed it with a larger data set up to $\itrpt{S}{k}{}$ to predict from a higher $t_k$), the PINN is \textit{re-trained}, and the resulting weights might slightly differ from those used to predict from a different IC (lower $t_k$).

In Figures~\ref{fig:rand-l2-error-interpolation} and \ref{fig:rand-l2-error-forecast}, the interpolation errors and forecast error graphs are shown for the \textbf{RAND} data sets, without any modifications or fine-tuning to either of the compared methods. It is clear that the PINN, which had previously performed very well on the \textbf{GRID} data sets, essentially failed and produced large $L_2$ relative errors. While it can be argued that the PINNs can be modified with respect to architecture, training schemes, activation functions, etc., this only highlights a more profound problem (and an open secret) for PINNs: their practical performance is highly sensitive to a multitude of hyperparameters and requires empirical ``rules of thumb'' (or worse, a lot of ``trial and error'') to make them work for specific cases or data layouts. On the other hand, our proposed DCBR is much more robust to changes in the training data layout, without any change in its settings. In fact, this is also unsurprising because of the principles-first manner in which we built DCBR.

Table~\ref{table:compute-times} summarizes the computation times of PINN and DCBR for the interpolation task. As expected, PINN requires substantial training time, even with CUDA acceleration (and is considerably slower on a CPU). In contrast, DCBR uses lazy evaluation, meaning that all computations are performed at query time. While each query is slower than evaluating a trained PINN, the cost remains modest (on a single-threaded CPU), and without any upfront (training) costs. This makes the effective time to the first solution much lower, which is suitable for real-time use in some settings, such as fluid dynamics. Moreover, DCBR does not need to be retrained when the data set is modified or enlarged, which enables online adaptation.

\begin{table}[tb]
\caption{Computational times (interpolation): \\pre-training, and per-query call}
{\renewcommand{\arraystretch}{1.2}
\begin{tabular}{|l|ll|l|}
\hline
\multirow{2}{*}{\textbf{Time}} & \multicolumn{2}{l|}{\textbf{PINN}}         & \multirow{2}{*}{\textbf{DCBR}} \\ \cline{2-3}
                      & \multicolumn{1}{l|}{CPU}  & CUDA  &                       \\ \hline
Pre-training (s)      & \multicolumn{1}{l|}{2056.16}     & 71.16 & 0.0                   \\ \hline
Per-query (ms)        & \multicolumn{1}{l|}{0.01} & 0.01  & 17.84                 \\ \hline
\end{tabular}
}
\label{table:compute-times}
\end{table}

\section{Conclusions}
\label{sec:conclusions}

We have introduced a new idea on how physics-based information can be embedded in data-driven regression by framing the prediction as an optimisation problem, enforcing the consistency of Taylor series expansions among samples via linear inequality constraints, and the PDE via equality constraints. Our case study on a reaction-diffusion system has shown competitive prediction performance with respect to physics-informed neural networks but without the need for a long training process. Furthermore, we have demonstrated that DCBR performance is robust even when presented with a different training data layout, whilst a physics-informed neural network has become unusable in terms of accuracy. Interesting further work on this topic could include the treatment of noisy measurements and uncertainty quantification, tailored numerical techniques to speed up computations, testing on more complex benchmark functions, and applications to realistic spatio-temporal systems.

\begin{ack}
The first author thanks Roberto Boffadossi from Politecnico di Milano for the insightful discussions that led to the seed idea for this paper.
\end{ack}

\section*{DECLARATION OF GENERATIVE AI AND AI-ASSISTED TECHNOLOGIES IN \\THE WRITING PROCESS}
During the preparation of this work the author(s) used ChatGPT 5.1 for information gathering and proofreading. After using this tool/service, the author(s) reviewed and edited the content as needed and take(s) full responsibility for the content of the publication.

\bibliography{ifacconf}             
                                                   







\end{document}